\documentclass[a4paper, 11pt,sumlimits,intlimits]{amsart}
\usepackage{amsmath,amsfonts,amsthm,amssymb,mathrsfs,latexsym,paralist, xypic, stmaryrd, soul}

\usepackage{tikz}
\usepackage{hyperref}

\usetikzlibrary{arrows,automata}

\tikzstyle{V}=[fill=black,circle,scale=0.4, outer sep = 4pt]

\allowdisplaybreaks

\newtheorem{thm}{Theorem}[section]

\theoremstyle{remark}

\newtheorem{example}[thm]{Example}
\newtheorem{examples}[thm]{Examples}
\theoremstyle{definition}
\newtheorem{defn}[thm]{Definition}

\newcommand{\bi}{\begin{itemize}}
\newcommand{\ei}{\end{itemize}}
\newcommand{\be}{\begin{enumerate}}
\newcommand{\ee}{\end{enumerate}}

\newcommand{\K}{\mathcal{K}}

\newcommand{\Z}{\mathbb{Z}}

\newcommand{\cstaralg}{$C^*$-algebra}

\begin{document}
\title{The UCT problem for nuclear $C^\ast$-algebras}
		
	\author{Nathanial Brown}\email{npb2@psu.edu}
		
	\author{Sarah L. Browne}\email{slbrowne@ku.edu}
	
	\author{Rufus Willett}\email{rufus@math.hawaii.edu}
	
	\author{Jianchao Wu}\email{jwu@math.tamu.edu}

\thanks{NB, SB and JW were partially supported by NSF grant \#DMS--1564401. RW was partially supported by NSF grant \#DMS--1564281 and \#DMS--1901522.}

\date{\today}
\maketitle
\begin{abstract}
In recent years, a large class of nuclear $C^\ast$-algebras have been classified, modulo an assumption on the Universal Coefficient Theorem (UCT). We think this assumption is redundant and propose a strategy for proving it. Indeed, following the original proof of the classification theorem, we propose bridging the gap between reduction theorems and examples. While many such bridges are possible,   various approximate ideal structures appear quite promising. 
\end{abstract}

\tableofcontents

\section{Introduction}

After decades of work by many hands, a remarkable classification theorem for simple $C^\ast$-algebras has emerged.\footnote{All $C^\ast$-algebras in this note are assumed separable, with exceptions like multiplier algebras.} Specifically, assuming the Universal Coefficient Theorem (UCT), when two simple  $C^*$-algebras have finite nuclear dimension, they are isomorphic if and only if they have isomorphic $K$-theoretic invariants (\cite{GongLinNiu-1,GongLinNiu-2}, \cite{Elliott:2015fb}, \cite{Tikuisis:2015kx}, \cite{Carrion:2020aa}); without finite nuclear dimension, classification via these invariants is impossible (\cite{Rordam:2003aa}, \cite{Toms:2008aa}). Thus the classification of simple nuclear $C^\ast$-algebras is complete -- modulo the UCT. This stunning fact has renewed interest in the old problem of whether or not every nuclear $C^\ast$-algebra satisfies the UCT. The purpose of this note is to review what is known about the UCT and propose a strategy for proving it for all nuclear $C^\ast$-algebras. 

The UCT is topological in nature, having its roots in Kasparov's $KK$-theory.  Kasparov introduced the $KK$-group $KK(A,B)$ around 1980 ~\cite{Kas80} for the purpose of building and analyzing maps between the $K$-theory groups of $C^*$-algebras $A$ and $B$.  In the Cuntz picture~\cite{Cun83}, an element of $KK(A,B)$ is represented by a quasi-homomorphism $A \rightarrow B$ and hence gives rise to a morphism at the level of $K$-theory $K_\ast(A) \rightarrow K_\ast(B)$. This induces a group homomorphism \[ \gamma \colon KK_\ast(A,B) \rightarrow \text{Hom}(K_\ast(A), K_\ast(B)).\] Unfortunately $\gamma$ cannot be an isomorphism, in general, since the left and right hand sides treat short exact sequences differently. Determining if $\gamma$ is surjective or describing its kernel is the role of the UCT.

In the stable case and Ext picture (\cite{Kas80}), an element of $KK_1(A,B)$ is a $C^\ast$-algebraic short exact sequence $0 \to B \to E \to A \to 0$. The boundary maps in the six-term exact sequence for $K$-theory then provide the homomorphism $K_\ast(A) \rightarrow K_{\ast +1}(B)$. When these maps vanish, i.e., belong to the kernel of $\gamma$, the $K$-theory of $E$ provides an element of $\text{Ext}(K_{\ast}(A), K_{\ast + 1}(B))$. Following the seminal paper of Rosenberg and Schochet \cite{RS87}, we define the {\em UCT class} to be those $C^\ast$-algebras $A$ for which \[ 0 \rightarrow \text{Ext}(K_{\ast}(A), K_{\ast + 1}(B)) \rightarrow KK_\ast(A,B) \stackrel{\gamma}{\rightarrow} \text{Hom}(K_\ast(A), K_\ast(B)) \rightarrow 0 \] is exact 
for any $C^*$-algebra $B$.\footnote{The reader is warned that our discussion only conveys broad ideas, and sweeps substantial and subtle details under the rug. Please see \cite{RS87} for a precise treatment.} $C^*$-algebras in the UCT class are said to {\em satisfy the UCT.}

The strategy we propose for proving nuclear $C^*$-algebras satisfy the UCT is rooted in the original proof of the stably finite case of the classification theorem for simple $C^*$-algebras of finite nuclear dimension. For a couple of decades, classification was only achieved for $C^*$-algebras constructed as limits of well-understood building blocks. It started with pioneering work of Elliott \cite{Elliott1976} on the classification of AF-algebras, i.e., limits of direct sums of matrix algebras. With increasingly sophisticated techniques, classification was achieved for simple A$\mathbb{T}$-algebras of real rank zero, and then, unital simple AH-algebras with no dimension growth (\cite{EG1996, Gong2002, EGL2007}). In 2000, Huaxin Lin made a conceptual leap with the introduction and classification of 
TAF-algebras (\cite{Lin01}), which are defined by an abstract approximation property as opposed to concrete inductive limit structure. This breakthrough was eventually generalized, leading to the classification of algebras with generalized tracial rank (g-TR) at most one (\cite{GongLinNiu-1,GongLinNiu-2}). Thus, the classes of examples which could be classified via $K$-theoretic invariants grew over time, getting larger with each decade. 

At the same time, very general approximation properties were being introduced and studied (\cite{Winter:2003kx}, \cite{Kirchberg:2004uq}, \cite{Winter:2010eb}), which led to reduction theorems in the classification program. That is, it was shown that in order to classify algebras in a large class, it suffices to classify a smaller subclass. Perhaps the most influential reduction theorem was due to Winter, who proved that to classify all simple $C^*$-algebras $A$ with finite nuclear dimension, it suffices to classify $A\otimes {\mathcal U}$ where ${\mathcal U}$ is the universal UHF algebra (\cite{Winter:2014aa}). This reduction is quite surprising as algebras of the form $A\otimes {\mathcal U}$ have several special properties not enjoyed in the general finite-nuclear-dimension case such as an abundance of projections and divisible $K$-theory. In the presence of other conditions like real rank zero and quasidiagonality, related reduction theorems inched closer and closer to the abstract approximation properties being classified (\cite{Winter:2006aa}).  In 2015, a remarkable bridge was constructed by Elliott, Gong, Lin and Niu (\cite{Elliott:2015fb}): if $A$ is simple, unital, satisfies the UCT, has finite nuclear dimension and every tracial state on $A$ is quasidiagonal, then $A\otimes {\mathcal U}$ has generalized tracial rank at most one.

We think following this roadmap could lead to a proof of the UCT for all nuclear C$^*$-algebras.\footnote{During the preparation of this article, Huaxin Lin proposed a different strategy in conference lectures. His interesting ideas will not be covered here.} In the next section we review existing reduction theorems. In section \ref{sec:examples} we review the classes of examples known to satisfy the UCT. In section \ref{sec:bridges} we discuss possible bridges between reduction theorems and examples: 
we focus on a new reduction theorem (Theorem~\ref{uct main}; also see Theorem~\ref{uct bridge}) using a class of $C^*$-algebras that we call \emph{decomposable} and that is related to the property of having finite nuclear dimension (see Example~\ref{wdec ex}). 
Finally, in section \ref{sec:dec uct} we sketch a proof of the UCT for decomposable $C^*$-algebras.

\section{Reduction theorems}\label{sec:reduction}

The first important subclass to which the UCT can be reduced is the so-called \emph{Kirchberg algebras}, i.e., simple nuclear and \emph{purely infinite} $C^*$-algebras. Recall that a simple $C^*$-algebra is purely infinite if every nonzero hereditary sub-{\cstaralg} contains an infinite projection, that is, a projection which is Murray-von Neumann equivalent to a proper subprojection of itself. Kirchberg algebras enjoy many useful properties, including:
\begin{enumerate}
	\item \label{enum-Kirchberg-properties-comparison-other} any two nonzero positive elements are Cuntz equivalent  (\cite[Proposition~4.1.1]{RS02}); 
	\item \label{enum-Kirchberg-properties-real-rank} real rank zero, meaning any self-adjoint element can be approximated in norm by self-adjoint elements with finite spectrum (\cite[Proposition~4.1.1]{RS02}); 
	\item \label{enum-Kirchberg-properties-O-infty} tensorial absorption of the Cuntz algebra $\mathcal{O}_\infty$ and thus the Jiang-Su algebra $\mathcal{Z}$ (\cite[Theorem~7.2.6]{RS02}), 
\end{enumerate}

Using the fact that $A$ is $KK$-equivalent to $A \otimes \mathcal{O}_\infty$ and Kirchberg's celebrated $\mathcal{O}_2$-embedding theorem (\cite[Theorem 6.3.11]{RS02}), an inductive limit construction shows that every nuclear C$^*$-algebra is $KK$-equivalent to a Kirchberg algebra (\cite[Proposition 8.4.5]{RS02}). Hence we have: 

\begin{thm}[Kirchberg]\label{thm:Kirchberg}
The UCT holds for all nuclear $C^*$-algebras if and only if it holds for all unital Kirchberg $C^*$-algebras.  
\end{thm}

In fact, Kirchberg reduced even further, to the simplest possible $K$-groups.  

\begin{thm}(\cite[Corollary 8.4.6]{RS02})\label{kirchberg k=0}
The UCT holds for all nuclear $C^*$-algebras if and only if it holds for all unital Kirchberg $C^*$-algebras with trivial $K$-theory (i.e., they are all isomorphic to $\mathcal{O}_2$).
\end{thm}

Our second reduction theorem deals with algebras at the opposite end of the spectrum from those which are simple and purely infinite. 

\begin{defn}
A {\cstaralg} is called \emph{RFD} or \emph{residually finite-dimensional} if it embeds into $\prod M_{k_i}$ for a sequence $(k_i)_{i \in \mathbb{N}}$ of integers, where $M_{k_i}$ denotes the {\cstaralg} of $k_i \times k_i$-matrices.
\end{defn}

Equivalently, a {\cstaralg} is RFD if it has a separating family of finite-dimensional representations. 
Using Voiculescu's stunning result that cones are always quasidiagonal (\cite{Voiculescu:1991aa}), Dadarlat established the following reduction theorem. 

\begin{thm}(\cite[proof of Lemma 2.4]{Dadarlat2003})\label{thm:dadarlat}
The UCT holds for all nuclear $C^*$-algebras if and only if it holds for all nuclear RFD $C^*$-algebras. 
\end{thm}

For our third reduction theorem we need Lin's groundbreaking tracial-approximation idea. 

\begin{defn}(\cite[page 694]{Lin01})\label{defn:TAF}
	A {\cstaralg} is \emph{TAF} or \emph{tracially approximately finite-dimensional} if for any $\epsilon > 0$, any finite subset $\mathcal{F} \subset A$ containing a non-zero element, and any full $a\in A_+$, there exists a finite dimensional $C^*$-subalgebra $B \subset A$ with $1_{B} = p$ such that, for all $x \in \mathcal{F}$, we have 
	\begin{enumerate}
		\item $\left\| p x - x p \right\| < \epsilon$, 
		\item the distance from $pxp$ to $B$ is no more than $\epsilon$, and
		\item $n [1-p] \leq [p]$ in the Murray-von Neumann semigroup of $A$ and $1-p$ is equivalent to a projection in the hereditary subalgebra generated by $a$. 
	\end{enumerate}
\end{defn}

The following theorem of Dadarlat, which relies on Theorem \ref{thm:dadarlat} and utilizes another inductive limit construction, is an analogue of Theorem~\ref{thm:Kirchberg} in the stably finite case. 

\begin{thm}(\cite[Theorem 1.2]{Dadarlat2003})
The UCT holds for all nuclear $C^*$-algebras if and only if it holds for all nuclear TAF $C^*$-algebras. 
\end{thm} 

In fact (\cite[Theorem 1.2]{Dadarlat2003}), analogously to Theorem \ref{kirchberg k=0}, it suffices to prove the UCT for a simple, unital, nuclear TAF $C^*$-algebra with the same $K$-theory as the universal UHF algebra $\mathcal{Q}$ (and so that it is isomorphic to $\mathcal{Q}$).

Our last reduction theorem requires another groundbreaking idea: noncommutative topological covering dimension. 

\begin{defn}(\cite[Definition 2.1]{Winter:2010eb})\label{defn:dimnuc}
	The \emph{nuclear dimension} of a separable {\cstaralg} $A$ is the infimum of all natural numbers $d$ such that there is a sequence $(F_i)_{i \in \mathbb{N}}$ of finite-dimensional {\cstaralg}s, a sequence of completely positive contractions $(\psi_i \colon A \to F_i)_{i \in \mathbb{N}}$, and $(d+1)$ sequences $(\phi^{(l)}_i \colon F_i \to A)_{i \in \mathbb{N}}$ of order-zero completely positive contractions for $l = 0,1, \ldots, d$ such that 
	\[
	\left\| \left( \phi^{(0)}_i + \ldots + \phi^{(d)}_i \right) \circ \psi_i (a) - a \right\| \to 0 \ \text{ as } i \to \infty
	\] 
	for any $a \in A$. 
\end{defn}

Note that having finite nuclear dimension implies nuclearity.  Every Kirchberg algebra has nuclear dimension one (\cite[Theorem G]{BBSTWW19}), hence Theorem \ref{thm:Kirchberg} implies the following. 

\begin{thm}\label{thm:dimnuc-1}
The UCT holds for all nuclear $C^*$-algebras if and only if it holds for all simple unital $C^*$-algebras with nuclear dimension 1. 
\end{thm}

There are several other open problems which are equivalent to the UCT for nuclear C$^*$-algebras. Since they are not reduction theorems in the sense we are considering the reader is referred to (\cite[Introduction]{Barlak:2017ac}) for a nice summary. 

\section{Examples}\label{sec:examples}

Rosenberg and Schochet observed that abelian C$^*$-algebras satisfy the UCT in (\cite[page 439]{RS87}). Essentially every other known example is derived from this case using a variety of permanence properties. In this section we recall the main examples, then review the long list of permanence properties enjoyed by the UCT class. 

A C$^*$-algebra is \emph{type I} if its double dual is a type I von Neumann algebra. Basic examples include abelian C$^*$-algebras and the compact operators on a Hilbert space. Another important class of examples, particularly for classification, are subhomogeneous C$^*$-algebras, i.e., subalgebras of $C(X) \otimes M_n({\mathbb C})$ for some space $X$ and $n \in \mathbb{N}$. 

\begin{thm} (\cite[page 439]{RS87})\label{thm:typeI} 
Type I C$^*$-algebras satisfy the UCT. 
\end{thm} 

Groupoid\footnote{All groupoids here are locally compact, Hausdorff, and second countable.} $C^*$-algebras (\cite{Renault:1980fk}, \cite{Williams:2019aa}) provide our next class of examples. For simplicity, we will stick to the amenable case (as it corresponds to nuclearity, at least for \'{e}tale groupoids). 

\begin{thm}(\cite{Tu99})\label{thm:groupoid}
Let $G$ be an amenable groupoid. Then its $C^\ast$-algebra satisfies the UCT. 
\end{thm}

Importantly, the previous result was partially generalized by Barlak and Li (\cite{BL17}) to \emph{twisted} \'{e}tale groupoid C$^*$-algebras. This allowed a connection with the notion of Cartan subalgebras. 

\begin{defn}\label{defn:Cartan} A maximal abelian self-adjoint subalgebra $B \subset A$ is called a \emph{Cartan} subalgebra if its normalizer generates $A$, it is the image of a conditional expectation, and it contains an approximate unit for $A$.
\end{defn} 

By work of Renault (\cite{Renault:2008if}), Cartan subalgebras induce twisted groupoid structures (which are amenable in the nuclear case) and hence we have the following: 

\begin{thm}\label{thm:Cartan} 
If $A$ is nuclear and has a Cartan subalgebra, then $A$ satisfies the UCT. 
\end{thm} 

The converse is also true when $A$ is simple and has finite nuclear dimension (\cite{Spielberg:2007aa}, \cite{Li:2018yv}).

Though narrow in scope when compared to the previous examples, Eckhardt and Gillaspy used special properties of nilpotent groups to prove the following interesting theorem (\cite{EG16}).

\begin{thm}\label{thm:nilpotent}
Let $G$ be a finitely generated nilpotent group and $\pi$ an irreducible representation of $G$. Then $C^\ast_\pi(G)$ satisfies the UCT.
\end{thm}
 
\subsection{Permanence Properties} 

Here are the known permanence properties of the nuclear UCT class. Definitions and descriptions of their utility follow. 
\begin{itemize} 
\item $KK$-equivalence 
\item Tensor products 
\item Inductive limits
\item Two out of three in a short exact sequence
\item Crossed products by $\mathbb{Z}$ or $\mathbb{R}$
\item Internal approximation by subalgebras
\end{itemize} 

C$^*$-algebras $A$ and $B$ are said to be $KK$-equivalent if there exists an invertible element in $KK(A,B)$. These equivalence classes are large, giving one lots of room to explore when searching for new examples within a class. For instance, the Kirchberg algebra $\mathcal{O}_\infty$ is $KK$-equivalent to $\mathbb{C}$! More generally, Rosenberg and Schochet proved that a $C^*$-algebra is in the UCT class if and only if it is $KK$-equivalent to an abelian C$^*$-algebra (\cite[pages 455-456]{RS02}, and also \cite[Proposition 5.3]{Skandalis:1988rr}). 

When $A$ and $B$ satisfy the UCT, so does their minimal (and maximal) tensor product $A \otimes B$, since $A$ and $B$ are $KK$-equivalent to abelian C$^*$-algebras. In particular, the stabilization of something in the UCT class remains in the UCT class. 

If $A_1 \to A_2 \to A_3 \to \cdots$ is an inductive system and each $A_i$ satisfies the UCT, then so does their inductive limit (\cite[Proposition 2.3]{RS87}). It follows that AF algebras, and their generalizations using subhomogeneous building blocks, satisfy the UCT. 

If $0 \rightarrow A \rightarrow D \rightarrow B \rightarrow 0$ is short exact and two of the algebras $A,D$ or $B$ satisfy the UCT, so does the third. In particular, the UCT class is closed under extensions and taking quotients by ideals in the UCT class (\cite[Proposition 2.3]{RS87}).

If $A$ satisfies the UCT and $\alpha$ is an action of either $\mathbb{Z}$ or $\mathbb{R}$, then the crossed products $A \rtimes_{\alpha} \mathbb{Z}$ or $A \rtimes_{\alpha} \mathbb{R}$ satisfy the UCT (\cite[Propositions 2.6 and 2.7]{RS87}). One can show Cuntz algebras satisfy the UCT this way, since their stabilizations are isomorphic to crossed products of AF algebras (\cite[page 87]{RS02}). 

The internal-approximation permanence property (which generalizes the inductive limit result) is a theorem of Dadarlat \cite[Theorem 1.1]{Dadarlat2003}. 

\begin{thm}
Let $A$ be a nuclear $C^\ast$-algebra. Assume for any finite set $\mathcal{F} \subset A$ and any $\epsilon >0$ there is a $C^\ast$-subalgebra $B$ of $A$ satisfying the UCT and such that $\text{dist}(a,B)<\epsilon$ for all $a \in \mathcal{F}$. Then $A$ satisfies the UCT.
\end{thm}

Taken together, these permanence properties are wide ranging and exceedingly useful. For instance, Tu's proof of the UCT for C$^*$-algebras associated to amenable groupoids first uses Kasparov's so-called Dirac-dual Dirac method to construct a $C^*$-algebra $A(G)$ which is $KK$-equivalent to $C^*(G)$. He then observes that $A(G)$ is an inductive limit of type I C$^*$-algebras, completing the proof. 

\section{Possible Bridges}\label{sec:bridges} 

Summarizing section \ref{sec:reduction}, we know the UCT holds for all nuclear C$^*$-algebras if and only if it holds for any of the following subclasses:   
\begin{itemize}
	\item Kirchberg algebras (with trivial K-theory);
	\item nuclear RFD algebras;
	\item simple, nuclear, unital TAF algebras;
	\item simple, unital C$^*$-algebras with nuclear dimension one.
\end{itemize} 

In section \ref{sec:examples} we saw that the following examples, and anything built out of them via appropriate permanence properties, satisfy the UCT:
\begin{itemize}
\item Type I C$^*$-algebras; 
\item $C^\ast(G)$, where $G$ is an amenable groupoid;
\item any C$^*$-algebra with a Cartan subalgebra. 
\end{itemize}

Any $KK$-equivalence from the first group to the second would prove the UCT for all nuclear C$^*$-algebras. For instance, one could try to prove that every Kirchberg algebra is $KK$-equivalent to something with a Cartan subalgebra. Or perhaps there is a notion of ``tracial Cartan subalgebra" which still allows one to prove the UCT, thereby adding another bullet point to the second group, and for which every TAF algebra is $KK$-equivalent to an algebra with this property. There are lots of possibilities. 

\subsection{Decomposable $C^*$-algebras}

In the classification program, bridging the gap between reduction theorems and examples took decades of hard work and experimentation. The same could be true for the UCT, but there is a potential bridge that seems particularly promising. It is based on ``approximate ideal structures'' as introduced in \cite{Willett:2019aa}, although in the very special case when the ``approximate ideals'' are finite dimensional.

\begin{defn}\label{decomp}
We say a unital $C^*$-algebra $A$ is \emph{decomposable}\footnote{It is shown in \cite{Jaime:2021vh} that the definition we give here is equivalent to the one from  \cite{Willett:2021te}.}. if for every $\epsilon > 0$ and every finite subset $X \subset A$ there is a triple $(h,C,D)$ consisting of a positive contraction $h$ in $A$ and finite-dimensional sub-$C^*$-algebras $C$ and $D$ of $A$ such that, under the operator norm and its induced metric, we have:
\begin{enumerate}
\item $\|[h,x]\|\leq \epsilon$ for all $x\in X$;
\item $d(hx,C)\leq \epsilon$ and $d((1-h)x,D)\leq \epsilon $ for all $x\in X$;
\item $d((1-h)hx,C\cap D)\leq \epsilon$ for all $x\in X$
\end{enumerate}
If one can only arrange the first two conditions, we say $A$ is a \emph{weakly decomposable} $C^*$-algebra.
\end{defn}

Observe that if $I, J \subset A$ are ideals such that $A= I+J$, an exercise using quasi-central approximate units shows that for any $\epsilon>0$ and finite $X\subset A$ there is a positive contraction $h\in A$ satisfying (1), (2), and (3) above with $I$ and $J$ in place of $C$ and $D$.  Philosophically then, $A$ is decomposable if it is ``approximately a sum of finite-dimensional ideals''. It is this analogy that gives applications to $KK$-theory via an approximate version of the Mayer-Vietoris exact sequence (see Section~\ref{sec:dec uct} and \cite{Willett:2021te}). On the other hand, despite this philosophical relation to ideals, there are many \emph{simple} examples.

\begin{examples}\label{dec ex}
The following simple $C^*$-algebras are decomposable:
\begin{itemize}
\item The Cuntz algebras $\mathcal{O}_n$ for $n$ finite (\cite{Jaime:2021vh}, following \cite[Section 7]{Winter:2010eb}). 
\item The crossed product $C(X)\rtimes \Z$ for any minimal action on a Cantor set (\cite[Section 2 and Section 8]{Guentner:2014aa}).
\end{itemize}
\end{examples}

It is not clear exactly how large the class of simple (nuclear) decomposable $C^*$-algebra is, but the above shows that it contains many interesting examples, and includes both infinite and finite simple $C^*$-algebras.  On the other hand, the class of simple $C^*$-algebras that are \emph{weakly} decomposable is certainly very large.  

\begin{example}[\cite{Jaime:2021vh}]\label{wdec ex}
Say $A$ has nuclear dimension one and real rank zero.  Then $A$ is weakly decomposable.\footnote{Conversely, if $A$ is weakly decomposable, then it has nuclear dimension one.  We do not know if (weak) decomposability implies real rank zero in general, but this is true in some interesting special cases such as when $A$ is separable and unital with unique trace \cite{Jaime:2021vh}.}
\end{example}

Combined with recent deep advances in the structure theory of simple nuclear $C^*$-algebras (\cite{Castillejos:2019ab}, \cite{Castillejos:2019aa}), this example implies that all simple, nuclear, $\mathcal{Z}$-stable $C^*$-algebras with real rank zero are weakly decomposable.   In particular, all Kirchberg algebras are weakly decomposable.  From this and Theorem \ref{kirchberg k=0} of Kirchberg, we get another reduction theorem.

\begin{thm}\label{thm:ApproxIdealReduction} If the UCT holds for all simple, unital, separable, weakly decomposable $C^*$-algebras with $K_*(A)=0$, then it holds for all nuclear C$^*$-algebras. 
\end{thm}

On the other hand, the third author and Yu have proved the following in \cite{Willett:2021te}, based on machinery built in \cite{Willett:2020aa}.  The proof will be explained in the next section.

\begin{thm}\label{uct the}
The UCT holds for all unital, separable, decomposable $C^*$-algebras with $K_*(A)=0$.
\end{thm}

As a direct consequence of this, the fact that $\mathcal{O}_2$ is decomposable (see Example~\ref{dec ex}), and Theorem \ref{kirchberg k=0} of Kirchberg, we get the following structural reformulation of the UCT.

\begin{thm}\label{uct main}
The following are equivalent: 
\begin{enumerate}
\item All nuclear $C^*$-algebras satisfy the UCT.
\item Any Kirchberg algebra with trivial $K$-theory is decomposable.
\end{enumerate}
\end{thm}

We would guess that not all Kirchberg algebras are decomposable (even though they are all weakly decomposable by Example \ref{wdec ex}).  Indeed, Jaime and the third author (\cite{Jaime:2021vh}) have shown that $C^*$-algebras satisfying a slight strengthening of decomposability have torsion-free $K_1$ group, and Kirchberg algebras can have arbitrary (countable) abelian groups as their $K$-theory.  We conjecture that a (UCT) Kirchberg algebra is decomposable if and only if it has torsion free $K_1$ group.

This section can be summarized by giving the following bridge to the UCT.  

\begin{thm}\label{uct bridge}
To prove the UCT for all nuclear $C^*$-algebras, it suffices to prove that any weakly decomposable $C^*$-algebra with trivial $K$-theory is equal to\footnote{or just $KK$-equivalent to} to a decomposable $C^*$-algebra.
\end{thm}

\section{The UCT for decomposable $C^*$-algebras}\label{sec:dec uct}

In this section, we will sketch the proof of Theorem \ref{uct the}.  The proof is technical; we give a different exposition here than in the original papers to try to make it more palatable.

We want to prove the UCT for a separable, unital, decomposable $C^*$-algebra $A$ with $K_*(A)=0$.  This is equivalent to showing that $KK(A,A)=0$.  Using a result of Dadarlat (\cite[Corollary 5.3]{Dadarlat:2005aa}), to show that $KK(A,A)=0$, it suffices to show that a certain quotient group $KL(A,A)$ of $KK(A,A)$ is $0$.  We will just use $KL$ as a black box in this exposition, so do not define it: suffice to say that it was introduced by R\o{}rdam (\cite[Section 5]{Rordam:1995aa}) in the case that $A$ satisfies the UCT, and in general it is defined as the largest Hausdorff quotient group of $KK(A,A)$ for a certain canonical topology on $KK(A,A)$ (\cite[Section 5]{Dadarlat:2005aa}).

The first step (\cite{Willett:2020aa}) in the proof of Theorem \ref{uct the} is to produce a good model for $KL(A,B)$, with $A$ and $B$ separable, unital, nuclear $C^*$-algebras.   To explain this, we first recall a non-standard description of the $K$-theory of $B$.  Let $B\otimes \K$ denote the stabilization of $B$, and let $\mathcal{M}(B\otimes \K)$ be its multiplier algebra.  Let $\mathcal{P}(B)$ be the set of all projections $p\in M_2(\mathcal{M}(B\otimes \K))$ such that $p-\begin{pmatrix} 1 & 0 \\ 0 & 0 \end{pmatrix}$ is in $M_2(B\otimes \K)$.  Then $K_0(B)$ canonically identifies with the set $\pi_0(\mathcal{P}(B))$ of path components of $\mathcal{P}(B)$.

Let now $\pi:A\to \mathcal{B}(\ell^2)$ be a unital, faithful, infinite multiplicity\footnote{This just means we take some faithful unital representation and add it to itself infinitely many times.} representation.  This gives rise to an inclusion 
$$
A\to \mathcal{B}(\ell^2)=\mathcal{M}(\K)\subseteq \mathcal{M}(B\otimes \K),
$$
which we use to identify $A$ with a $C^*$-subalgebra of $\mathcal{M}(B\otimes \K)$ (and therefore with a subalgebra of $M_2(\mathcal{M}(B\otimes \K))$ by having it act diagonally).  For a subset $X\subset A$ and $\epsilon\geq 0$, define  
\begin{equation}\label{pespb}
\mathcal{P}_\epsilon(X,B):=\{p\in \mathcal{P}(B)\mid \|[p,x]\|\leq \epsilon \text{ for all } x\in X\}.
\end{equation}
Let $KK_\epsilon(X,B):=\pi_0(\mathcal{P}_\epsilon(X,B))$ be the set of path components of $\mathcal{P}_\epsilon(X,B)$, which is an abelian group in a natural way.  As discussed above, $K_0(B)=\pi_0(\mathcal{P}_\epsilon(\varnothing,B))$ for any $\epsilon$, so one can think of $KK_\epsilon(X,B)$ as ``the part of $K_0(B)$ that commutes with $X$ up to $\epsilon$ error''.

Let now $(X_n)$ be an increasing sequence of finite subsets of $A$ with dense union and let $(\epsilon_n)$ be a decreasing sequence of positive numbers tending to zero.  Then there is a system of homomorphisms
\begin{equation}\label{inv sys}
\cdots \to KK_{\epsilon_n}(X_n,B)\to KK_{\epsilon_{n-1}}(X_{n-1},B) \to \cdots \to KK_{\epsilon_1}(X_1,B),
\end{equation}
where each arrow comes from the fact that commuting with $X_{n-1}$ up to $\epsilon_{n-1}$ is easier than commuting with $X_{n}$ up to $\epsilon_{n}$.   Hence we can define the inverse limit ${\displaystyle \varprojlim KK_{\epsilon_n}(X_n,B)}$ in the sense of abelian group theory.  One can also (we will not go into the details as we do not want to define $KL(A,B)$) construct homomorphisms $\kappa_n:KL(A,B)\to KK_{\epsilon_n}(X_n,B)$ for each $n$.

\begin{thm}[\cite{Willett:2020aa}]
With notations as above, the maps $\kappa_n$ fit together to define an isomorphism
$$
KL(A,B) \xrightarrow{\ \cong\ } \varprojlim KK_{\epsilon_n}(X_n,B).
$$
\end{thm}

In particular, this theorem implies (by definition of the inverse limit) that to show that $KL(A,B)=0$, it suffices to show that for any $\alpha\in KL(A,B)$ and any $n$, there exists $N\geq n$ such that $\kappa_N(\alpha)$ maps to zero under the canonical map $KK_{\epsilon_N}(X_N,B)\to KK_{\epsilon_n}(X_n,B)$.\\

We now turn to the second step (\cite{Willett:2021te}) in the proof of Theorem \ref{uct the}, which is a Mayer-Vietoris argument.  Our aim is to show that for any class $\alpha\in KL(A,A)=0$ and any $n$, there exists $N\geq n$ such that $\kappa_N(\alpha)$ maps to zero in $KK_{\epsilon_n}(X_n,A)$.     Recall\footnote{We cannot find this in the literature, but it is well-known and can be derived from the usual long exact sequence by the arguments of \cite[Proposition 2.7.15]{Willett:2010ay}, for example.} that if $A$ is nuclear and splits as a sum of two ideals $A=I+J$, then there is a six-term exact Mayer-Vietoris sequence
$$
\cdots \to KK(I\cap J,SA) \stackrel{\partial}{\to} KK(A,A) \stackrel{\sigma}{\to} KK(I,A)\oplus KK(J,A) \to \cdots 
$$
where $SA$ is the suspension of $A$.  We want to reproduce this ``locally''.  

Choose a triple $(h,C,D)$ with the properties in Definition \ref{decomp} for the given $X_n$ and some suitably small $\epsilon$ (determined by the proof).  We claim that for any large enough $N$ (how large depends on $n$, $h$, $C$, $D$) we can find a collection of maps
$$
\xymatrix{ &  KK_{\epsilon_N}(X_N,A) \ar[r]^-\sigma \ar[d] & KK(C,A)\oplus KK(D,A) \\
KK(C\cap D,SA) \ar[r]^\partial & KK_{\epsilon_n}(X_n,A) & },
$$
where the vertical map comes from line \eqref{inv sys}, and with the following exactness property: if an element of $KK_{\epsilon_N}(X_N,A)$ is mapped to zero in $KK(C,A)\oplus KK(D,A)$ by $\sigma$, then its image in $KK_{\epsilon_n}(X_n,SA)$ comes from $KK(C\cap D,SA)$ via $\partial$.  However, the groups $KK(C,A)$, $KK(D,A)$ and $KK(C\cap D,SA)$ are all zero as $C$ and $D$ are finite dimensional and $K_*(A)=0$, so the claim completes the proof.

The construction of the ``local Mayer-Vietoris sequence'' in the claim above uses the same ideas as the classical six-term exact sequence in $K$-theory.   The map 
$$
\sigma: KK_{\epsilon_N}(X_N,A)\to KK(C,A)\oplus KK(D,A).
$$
is defined by choosing $X_N$ large enough to (approximately) contain the unit balls of $C$ and $D$.  A cycle for $KK_{\epsilon_N}(X_N,A)$ thus approximately commutes with the unit ball of $C$, and is therefore close to an element that actually commutes with $C$ by averaging over the (compact!) unitary group of $C$.  For $C$ finite-dimensional, one can show that the group $KK(C,A)$ identifies\footnote{This is cheating in general: there is a mild extra complication if $C$ is not a unital subalgebra of $A$.  However, it is precisely true when $K_*(A)=0$, which is all we need.} with the set of path components of $\mathcal{P}_0(C,A)$ (see line \eqref{pespb} for notation), so we get a map $KK_{\epsilon_N}(X_N,A)\to KK(C,A)$.  A map $KK_{\epsilon_N}(X_N,A)\to KK(D,A)$ is defined similarly, and putting these together defines $\sigma$.  

The map 
$$
\partial: KK(C\cap D,SA) \to KK_{\epsilon_n}(X_n,A)
$$
is defined by adapting classical, purely algebraic, formulas\footnote{These go back at least to Milnor's exposition of algebraic $K$-theory (\cite[Chapter 2]{Milnor:1971kl}).} for the boundary map in $K$-theory in terms of $h$ and $1-h$.  Given an element $\beta$ in the kernel of $\sigma$, these formulas also give an ansatz for an element of $KK(C\cap D,SA)$ that maps to the image of $\beta$ in $KK_{\epsilon_n}(X_n,A)$.  A careful approximation (and using the same averaging argument as above, this time based on the compactness of the unitary group of $C\cap D$) shows that this ansatz can be made to work, completing the proof.

\end{document}